\documentclass[12pt,twoside]{amsart}
\usepackage{amssymb}

\newtheorem{theorem}{Theorem}[section]
\newtheorem{lemma}[theorem]{Lemma}
\newtheorem{proposition}[theorem]{Proposition}
\newtheorem{corollary}[theorem]{Corollary}

\theoremstyle{definition}
\newtheorem{definition}[theorem]{Definition} 

\newtheorem{example}[theorem]{Example}

\newcommand{\skipit}[1]{{}}
\newcommand{\prfend}{\hbox to7pt{\hfil}
\par\vskip-\baselineskip\hbox to\hsize
{\hfil\vbox {\hrule width6pt height6pt}}\vskip\baselineskip}

\newcommand{\myarrow}[2]{\hbox to #1pt{\hfil$\to$\hfil}{\hskip-#1pt{\raise
10pt\hbox to#1pt{\hfil$\scriptscriptstyle #2$\hfil}}}}

\title{Bounds and asymptotic minimal growth\\for Gorenstein Hilbert functions}

\author[Juan Migliore]{Juan Migliore${}^*$}
\address{
Department of Mathematics, University of Notre Dame, Notre Dame, IN
46556, USA}
\email{Juan.C.Migliore.1@nd.edu}
\author[Uwe Nagel]{Uwe Nagel${}^{+\,\#}$}
\address{Department of Mathematics,
University of Kentucky, 715 Patterson Office Tower,
Lexington, KY 40506-0027, USA}
\email{uwenagel@ms.uky.edu}

\author[Fabrizio Zanello]{Fabrizio Zanello${}^\#$}
\address{Department of Mathematical Sciences, Michigan
Technological University, Houghton, MI 49931-1295, USA}
\email{zanello@math.kth.se}

\thanks{${}^*$ Part of the work for this paper was done while the first
author was sponsored by the National Security Agency under Grant
Number H98230-07-1-0036.\\
${}^+$ Part of the work for this paper was done while the second
author was sponsored by the National Security Agency under Grant
Number H98230-07-1-0065.\\
${}^{\# }$ The idea of this work originated from a discussion between the second and the third author, during a workshop at the Institute for Mathematics and Applications (IMA). They both thank the IMA for partial support.
}

\begin{document}

\begin{abstract}
We determine new bounds on the entries of Gorenstein Hilbert functions, both in any fixed codimension and asymptotically.

Our first main theorem is a lower bound for the degree $i+1$ entry of a Gorenstein $h$-vector, in terms of its entry in degree $i$. This result carries interesting applications concerning unimodality: indeed, an important consequence is that, given $r$ and $i$, all Gorenstein $h$-vectors of codimension $r$ and socle degree $e\geq e_0=e_0(r,i)$ (this function being explicitly computed) are unimodal up to degree $i+1$. This immediately gives a new proof of a theorem of Stanley that all Gorenstein $h$-vectors in codimension three are unimodal.

Our second main theorem is an asymptotic formula for the least value that the $i$-th entry of a Gorenstein $h$-vector may assume, in terms of codimension, $r$, and socle degree, $e$. This theorem broadly generalizes a recent result of ours, where we proved a conjecture of Stanley predicting that asymptotic value in the specific case $e=4$ and $i=2$, as well as a result of Kleinschmidt which concerned the logarithmic asymptotic behavior in degree $i=\left \lfloor \frac{e}{2}\right \rfloor $.
\end{abstract}
{\ }\\

\maketitle

\section{Introduction}

It has been observed by Bass \cite{bass} that Gorenstein algebras are ubiquitous throughout mathematics.  Despite this, it is often distressingly difficult to find ones with specific desired properties (e.g. in liaison theory, to find ``small'' Gorenstein subschemes containing a given subscheme of projective space).  A first step is to have some understanding of which Hilbert functions can occur.  These are completely classified in codimension $r\leq 3$ \cite{St} (see also the third author's \cite{Za3}; Macaulay first proved the result in the simpler case $r=2$ in \cite{Ma}), but a complete description is unknown if the codimension is at least $4$, in spite of a remarkable amount of work performed by several authors (see, e.g., \cite{BI,BL,Bo,IS,MNZ1,MNZ2,St}).

A first step is to have some idea of the extremal general ``shape'' of the Hilbert function as the codimension gets arbitrarily large.   The upper bound is, of course, the case of compressed Gorenstein algebras (see Emsalem-Iarrobino's \cite{EI}, which was the seminal work on compressed algebras, and also \cite{Ia1,FL,Za1,Za2}, which further developed and extended this theory), so the interesting question is to ask for a lower bound.  This has first been done by Stanley \cite{St2}.  In particular, he considered the special case where the socle degree is 4, and gave a precise conjecture for the asymptotic growth of the least value, $f(r)$, of the Hilbert function  in degree 2, in terms of the codimension, $r$.  Specifically, he conjectured that
\[
\lim_{r \rightarrow \infty} \frac{f(r)}{r^{2/3}} = 6^{2/3}.
\]

This conjecture appeared in the first edition  of \cite{St2}, in
1983.  Bounds were given by Stanley and by Kleinschmidt
\cite{kleinschmidt}, but the precise limit was only proved
(verifying Stanley's conjecture) in 2006 by the current authors.

The purpose of this paper is to give a very broad  generalization of
this result.  There are at least two initial questions that one can
ask concerning the general shape of the Hilbert function of a
Gorenstein algebra, and those will be answered in this paper.
First, if we know the entry of a Hilbert function in any given
degree, what is a ``good'' lower bound for the value it can assume
in the next degree?

Our answer to this question, Theorem \ref{lower bd},  carries very
interesting applications concerning unimodality: indeed, an
important consequence of our result is that, given $r$ and $i$, all
Gorenstein $h$-vectors of codimension $r$ and socle degree $e\geq
e_0=e_0(r,i)$ (this constant being explicitly computed) are unimodal
up to degree $i+1$.

In codimension $r\leq 3$, this result is powerful  enough to supply
a new proof of a celebrated theorem of Stanley that all Gorenstein
$h$-vectors are unimodal.

Second, one can ask for asymptotic bounds given only  in terms of
the codimension and the socle degree, as in Stanley's situation.  In
Theorem \ref{main result}, we will supply the least asymptotic value
that the $i$-th entry of a Gorenstein $h$-vector may assume, in
terms of  codimension $r$ and socle degree $e$. This generalizes the
recent result of ours mentioned above, where we solved a conjecture
of Stanley predicting that asymptotic value in the specific case
$e=4$ and $i=2$, as well as a result of Kleinschmidt
(\cite{kleinschmidt}, Theorem 1) for $i=\left \lfloor
\frac{e}{2}\right \rfloor $.

Our asymptotic result follows by combining our lower bounds and a
construction of suitable Gorenstein algebras. We illustrate this
with a specific example.
\medskip

\noindent {\bf Example.} Consider the Gorenstein $h$-vectors of the
form
\[
(1, 125, h_2, h_3,\ldots, 125, h_8 = 1).
\]
The proof of Theorem \ref{main result} guarantees the existence of 
a Gorenstein $h$-vector \linebreak $(1, 125, 95, 77, 71, \dots)$.  On the other hand, with this value of $h_1$, Theorem \ref{lower bd} provides $h_2 \geq 95$, which is thus sharp.  Taking $h_2 = 95$, it provides $h_3 \geq 77$, again sharp thanks to our explicit construction.  Then taking $h_3 = 77$, we obtain $h_4 \geq 70$, but now our example does not achieve the bound; indeed, we do not know if this bound is sharp or not.  See also Example \ref{relate results} for a broad generalization of this example.

On the other hand, it is not surprising that our bound in Theorem \ref{lower bd} is not always 
 sharp, since a sharp bound would probably make it easy to decide
if non-unimodal Gorenstein $h$-vectors of codimension four do exist.  See also Example \ref{ex-sharpness}.
What we do find surprising is that Theorem \ref{lower bd} is strong
enough to give a sharp asymptotic bound (Theorem \ref{main result}),
as described above.


\section{A general lower bound and its applications}

Throughout this paper, $k$ will denote an infinite field, and $R = k[x_1,\dots,x_r]$ a graded polynomial ring in $r$ variables.  Each standard graded $k$-algebra $A$ can be written as $A = R/I$, where  $I \subset R$ is  a homogeneous ideal.

We begin by recalling  results of Macaulay, Green, and Stanley that we will need in this paper.

\begin{definition} Let $n$ and $i$ be positive integers. The {\em i-binomial expansion of n} is
\[
n_{(i)} = \binom{n_i}{ i}+\binom{n_{i-1}}{i-1}+...+\binom{n_j}{j},
\]
 where $n_i>n_{i-1}>...>n_j \geq j \geq 1$. We remark that such an expansion always exists and it is unique (see, e.g.,  \cite{BH} Lemma 4.2.6).
\end{definition}

Following \cite{BG}, we define, for any integers $a$ and $b$,
\[
\left(n_{(i)}\right)_{a}^{b}=\binom{n_i+b}{i+a}+\binom{n_{i-1}+b}{i-1+a}+...+\binom{n_j+b}{j+a},
\]
where we set $\binom{m}{q}=0$ whenever $m<q$ or $q<0$.

\begin{theorem}\label{gr}  Let $L \in A$ be a general linear form. Denote by $h_d$  the degree $d$ entry of the Hilbert function of $A$ and by $h_d^{'}$ the degree $d$ entry of the Hilbert function of $A/L A$. Then:
\begin{itemize}
\item[i)] (Macaulay) $\displaystyle h_{d+1}\leq \left ( \left(h_d\right)_{(d)} \right )_{1}^{1}.$

\item[ii)]  (Green) $\displaystyle h_d^{'}\leq \left ( \left(h_d\right )_{(d)} \right )_{0}^{-1}.$
\end{itemize}
\end{theorem}

\begin{proof}
i) See \cite{BH}, Theorem 4.2.10.

ii) See \cite{Gr}, Theorem 1.
\end{proof}

The following  simple observation is not new (see for instance \cite{St2}, bottom of p. 67):

\begin{lemma}(Stanley)\label{stle}.  Let $A$ be an artinian Gorenstein algebra, and let $L\in A$ be any linear form. Then the  Hilbert function of $A$ can be written as
\[
h := (1,h_1,...,h_e) = (1, b_1+c_1,..., b_e+c_e = 1),
\]
 where
\[
b=(b_1=1,b_2,...,b_{e-1},b_e=1)
\]
is the $h$-vector of $A/(0:L)$ (with the indices shifted by 1), which is a Gorenstein algebra, and
\[
c=(c_0=1,c_1,...,c_{e-1},c_e=0)
\]
is the $h$-vector of $A/L A$.
\end{lemma}

Perhaps the most important, and definitely the most consequential, result of this section is a lower bound for the value of a Gorenstein Hilbert function in terms of the value in the previous degree. This result generalizes \cite{MNZ1}, Theorem 4, which treated the case $i=1$.  Specifically, we have that:

\begin{theorem} \label{lower bd}
Suppose that ${h} = (1,h_1 = r ,h_2,\dots,h_{e-2}, h_{e-1} = r, h_e = 1)$ is the $h$-vector of an artinian Gorenstein algebra over $R = k[x_1,\dots,x_r]$.  Assume that $i$ is an integer satisfying $1 \leq i \leq \frac{e}{2} -1$.  Then
\[
h_{i+1} \geq \left ( \left(h_i\right )_{(e-i)} \right )^{-1}_{-1} + \left ( \left ( h_i\right )_{(e-i)} \right )^{-(e-2i)}_{-(e-2i-1)}.
\]
\end{theorem}

\begin{proof} As in Stanley's Lemma \ref{stle}, let us write $h=(1,h_1,...,h_e)$ as $b+c=(1, b_1+c_1,..., b_e+c_e = 1)$, where we have picked the form $L$ to be general inside $R$. Notice that $b$ is a Gorenstein $h$-vector of socle degree $e-1$. Therefore, by symmetry and our choice of the indices, $b_j=b_{e+1-j}$ for all $j$.

Hence, by Green's Theorem \ref{gr}, ii), we have that $$c_{e-i}\leq \left( \left(h_{e-i}\right)_{(e-i)}\right)_0^{-1}=\left(\left(h_{i}\right)_{(e-i)}\right)_{0}^{-1}.$$

Thus (using the Pascal's Triangle inequality) $$b_{i+1}=b_{e-i}\geq h_i-\left(\left(h_i\right)_{(e-i)}\right)_{0}^{-1}=\left(\left(h_{i}\right)_{(e-i)}\right)_{-1}^{-1},$$
which implies 
$$c_{i+1}\leq h_{i+1}-\left(\left(h_{i}\right)_{(e-i)}\right)_{-1}^{-1}.$$

By iterating Macaulay's Theorem \ref{gr}, i), we obtain another upper bound on $c_{e-i}$, namely:
$$c_{e-i}\leq \left(\left(h_{i+1}-\left(\left(h_{i}\right)_{(i)}\right)_{-1}^{-1}\right)_{(i+1)}\right)^{e-2i-1}_{e-2i-1}.$$

Now, since by Green's theorem, $c_{e-i}\leq ((h_{i})_{(e-i)})_0^{-1}$, we write $c_{e-i}=((h_{i})_{(e-i)})_0^{-1}-a$, for some integer $a\geq 0$. Therefore,
$$c_{i+1}=h_{i+1}-b_{i+1}=h_{i+1}-b_{e-i}=h_{i+1}-h_i+\left(\left(h_{i}\right)_{(e-i)}\right)_0^{-1}-a=h_{i+1}-\left(\left(h_{i}\right)_{(e-i)}\right)_{-1}^{-1}-a,$$
from which we get, again by iterating Macaulay's theorem:
$$\left(\left(h_{i+1}-\left(\left(h_{i}\right)_{(e-i)}\right)_{-1}^{-1}-a\right)_{(i+1)}\right)^{e-2i-1}_{e-2i-1}=\left(\left(c_{i+1}\right)_{(i+1)}\right)^{e-2i-1}_{e-2i-1}\geq c_{e-i}=\left(\left(h_i\right)_{(e-i)}\right)_0^{-1}-a.$$
Therefore, since $(m_{(d)})^1_1$ is a strictly increasing function of $m$, we have
$$\left(\left(h_{i+1}-\left(\left(h_{i}\right)_{(e-i)}\right)_{-1}^{-1}\right)_{(i+1)}\right)^{e-2i-1}_{e-2i-1}\geq $$$$ \left(\left(h_{i+1}-\left(\left(h_{i}\right)_{(e-i)}\right)_{-1}^{-1}-a\right)_{(i+1)}\right)^{e-2i-1}_{e-2i-1}+a\geq \left(\left(h_{i}\right)_{(e-i)}\right)_0^{-1}.$$
In particular, we have proved that 
$$\left(\left(h_{i+1}-\left(\left(h_{i}\right)_{(e-i)}\right)_{-1}^{-1}\right)_{(i+1)}\right)^{e-2i-1}_{e-2i-1}\geq  \left(\left(h_{i}\right)_{(e-i)}\right)_0^{-1}.$$

From the last inequality, we again use monotonicity. Notice first that $((h_{i})_{(e-i)})_0^{-1}$
already presents itself in its $(e-i)$ binomial expansion (possibly after eliminating those binomial coefficients that are equal to 0). Similarly, the left-hand side of the last inequality is also already written as an $(e-i)$ binomial expansion, since $(i+1)+(e-2i-1)=e-i$. Hence we easily get
$$h_{i+1}-\left(\left(h_{i}\right)_{(e-i)}\right)_{-1}^{-1}\geq \left(\left(h_{i}\right)_{(e-i)}\right)_{-(e-2i-1)}^{-1-(e-2i-1)}=\left(\left(h_{i}\right)_{(e-i)}\right)_{-(e-2i-1)}^{-(e-2i)},$$
as we wanted to show.
\end{proof}

Let us now present some interesting applications of the above theorem.

\begin{example}
In \cite{BI}, D.\ Bernstein and A.\ Iarrobino gave the first known example of a non-unimodal Gorenstein Hilbert function in codimension 5, namely
\[
1 \ \ 5 \ \ 12 \ \ 22 \ \ 35 \ \ 51 \ \ 70 \ \ 91 \ \ 90 \ \ 91 \ \ 70 \ \ 51 \ \ 35 \ \ 22 \ \ 12 \ \ 5 \ \ 1.
\]
One can check, using Theorem \ref{lower bd}, that with  a value of
91 in degree 7, the 90 in degree 8 is optimal.  Notice that, for $i
\leq 6$, the value of $h_{i+1}$ is not optimal with respect to our
bound, but in any case Theorem \ref{lower bd} guarantees that the
Hilbert function is unimodal in those degrees; only in degrees 7 and
8 it might  be possible to violate unimodality.

This motivates Proposition \ref{unimodal} below.
\end{example}

\begin{proposition} \label{unimodal}
Let ${h} = (1,h_1,\dots , h_{e-1}, h_e)$ be a Gorenstein $h$-vector,
and fix an index $i \leq \frac{e}{2} -1$.  If
$$h_i<\binom{e-i+2}{2}-\binom{e-2i-1}{2}=\frac{1}{2}(i+3)(2e-3i),$$
then $h_{i+1} \geq h_i$.
\end{proposition}

\begin{proof}
Since $\binom{e-i+2}{2}=\binom{e-i+2}{e-i}$, the $(e-i)$-binomial expansion of $h_i$ clearly begins with $\binom{e-i+1}{e-i} + \dots$ or $\binom{e-i}{e-i} + \dots$.  In particular, since $h_i<\binom{e-i+2}{2}-\binom{e-2i-1}{2}$, there exists an integer $\ell ,$ $ e-2i-1 \leq \ell \leq e-i+1$, such that
\[
h_i = \binom{e-i+1}{e-i} + \dots + \binom{\ell +1}{\ell} + (*) = (e-i+1) + \dots + (\ell +1) + (*),
\]
where $(*)$ is a sum of at most $\ell -1$ binomial coefficients of the form $\binom{j}{j}$. (Of course, $h_i:=(*)$ if $\ell = e-i+1$.)  Theorem \ref{lower bd} then gives
\[
\begin{array}{rcl}
h_{i+1} & \geq & ((e-i) + \dots + \ell + (*)) + (e-i -\ell +1) \\
& = & (e-i+1) + (e-i) + \dots + (\ell +1) + (*) \\
& = & h_i,
\end{array}
\]
as desired.
\end{proof}

A simple, but very important, application of the previous result is the following:

\begin{corollary}\label{ccc}
Fix $r$ and $i$. Then all Gorenstein $h$-vectors of codimension $r$ and socle degree
$$e> \frac{(i+1)(i+2)\cdot \cdot \cdot (i+r-1)}{(i+3)(r-1)!}+\frac{3}{2}i$$
are unimodal up to degree $i+1$.
\end{corollary}
\begin{proof}
One can show that, for any $r\geq 2$, $$\frac{(i+1)(i+2)\cdot \cdot \cdot (i+r-1)}{(i+3)(r-1)!}+\frac{3}{2}i\geq \frac{i(i+1)\cdot \cdot \cdot (i+r-2)}{(i+2)(r-1)!}+\frac{3}{2}(i-1).$$

Hence it suffices to prove unimodality only in degree $i+1$.

Since $\frac{1}{2}(i+3)(2e-3i)$ is an increasing function of $e$, and $h_i$ clearly cannot exceed $\binom{i+r-1}{i}$, by Proposition \ref{unimodal} it is enough to show that $\frac{1}{2}(i+3)(2e-3i)>\binom{i+r-1}{i},$ and a standard computation shows that this is equivalent to the inequality on $e$ of the statement.
\end{proof}

In particular, our result is strong enough to reprove the well-known theorem of Stanley that all codimension $r\leq 3$ Gorenstein $h$-vectors are unimodal (see also \cite{Za3}):

\begin{theorem}[\cite{St}]\label{corst}
All Gorenstein $h$-vectors of codimension $r\leq 3$ are unimodal.
\end{theorem}

\begin{proof}
Notice that, a fortiori, it suffices to show that the inequality of the statement of Corollary \ref{ccc} is satisfied for $i=\left \lfloor \frac{e}{2}\right \rfloor -1$ and $r=3$. Therefore, we want to prove that
$$e>\frac{\left \lfloor \frac{e}{2}\right \rfloor \left(\left \lfloor \frac{e}{2}\right \rfloor +1\right)}{2\left(\left \lfloor \frac{e}{2}\right \rfloor +2\right)}+\frac{3}{2}\left(\left \lfloor \frac{e}{2}\right \rfloor -1\right).$$
But the right hand side is equal to $$2\left \lfloor \frac{e}{2}\right \rfloor -\frac {2\left \lfloor \frac{e}{2}\right \rfloor +3}{\left \lfloor \frac{e}{2}\right \rfloor +2 },$$ and  the desired inequality  immediately follows, since $e\geq 2\left \lfloor \frac{e}{2}\right \rfloor $.
\end{proof}

For $r=4$ the estimate we obtain is still a very interesting one. Namely, from Corollary \ref{ccc}, we immediately have:

\begin{corollary}\label{cor4}
All Gorenstein $h$-vectors of codimension $4$ and socle degree $e>\frac{1}{6}(i^2+12i+2)$ are unimodal up to degree $i+1$.
\end{corollary}

This complements the main result of \cite{MNZ2}, which  focused on
the initial degree of $I$ rather than on the socle degree of $R/I$.
There it was shown that, whenever $r=4$ and $h_4\leq 33$, then the
possible $h$-vectors for Gorenstein algebras are {\em precisely} the
SI-sequences.
\smallskip

We conclude this section with an example showing that the bound
given in Theorem \ref{lower bd} is not always sharp. However, in the
next section we will prove that this bound is asymptotically sharp.

\begin{example}
  \label{ex-sharpness}
Consider Gorenstein $h$-vectors of the form
\[
(1, 4, 10, 20, h_4, h_5, h_6 = h_4, 20, 10, 4, 1).
\]
Assume $h_4 = 33$. Then Theorem \ref{lower bd} gives $h_5 \geq 30$,
whereas Theorem 3.1 in \cite{MNZ2} says that $h_5 \geq h_4 = 33$. In
fact, using the methods of \cite{MNZ2}, Theorem 3.1, one can show
that all the above Gorenstein $h$-vectors are unimodal.

Notice that the methods in \cite{MNZ2} work nicely for algebras with
low initial degree whose codimension is at most four. The methods
developed in this paper work in general. This is the big advantage
of the current paper.
\end{example}


\section{Asymptotic minimal growth}

The following definition  generalizes one introduced in \cite{St2} and extended in \cite{MNZ1}.

\begin{definition} Fix integers $e$ and $i$.  Then $f_{e,i}(r)$ is the least possible value in degree $i$ of the Hilbert function of a Gorenstein algebra with socle degree $e$ and codimension $r$.
\end{definition}

\begin{lemma}[\cite{BG}, Lemma 3.3] \label{BG Lemma}
Let $A,d$ be positive integers.  Then

\begin{enumerate}
\item Assume that $d > 1$ and $s := \left (A_{(d)} \right )^{-1}_{-1}$.  Then $s$ is the smallest integer such that $A \leq \left ( s_{(d-1)} \right )^{+1}_{+1}$.

\item Assume that $d > i$.  Then
\[
\left ( A_{(d)} \right )^{-i-1}_{-i-1} = \left ( \left ( \left ( A_{(d)} \right )^{-i}_{-i} \right )_{(d-i)} \right )^{-1}_{-1}.
\]

\end{enumerate}
\end{lemma}

We note the following two immediate consequences of Lemma \ref{BG Lemma}.

\begin{corollary} \label{minus}
With the notation of Lemma \ref{BG Lemma} we have

\begin{enumerate}

\item $\displaystyle \left ( \left ( \left ( A_{(d)} \right )^{-1}_{-1} \right )_{(d-1)} \right )^{-1}_{-1} = \left ( A_{(d)} \right )^{-2}_{-2}. $

\item $\displaystyle  \left ( \left ( \dots\left (  \left ( \left ( \left ( A_{(d)} \right )^{-1}_{-1} \right )_{(d-1)} \right )^{-1}_{-1} \right ) \dots \right )_{(d-i)}  \right )^{-1}_{-1} = \left ( A_{(d)} \right )^{-(i+1)}_{-(i+1)}.$

\end{enumerate}
\end{corollary}

We need two more preliminary results before proving our main theorem. Remember that, given two functions $f$ and $g$, we say that $f(m)\in O(g(m))$ if, for $m$ large, there exists a positive constant $C$ such that $\vert f(m)\vert \leq C\cdot g(m)$.

\begin{lemma}\label{r}
Given $e\geq 1$, every positive integer $r$ can be written in the form
$$r=m+\binom{m+e-3}{e-1}+\binom{a_{e-2}}{e-2}+\binom{a_{e-3}}{e-3}+...+\binom{a_1}{1},$$
where $m$ is the largest integer such that $m+\binom{m+e-3}{e-1}\leq r$, $a_{e-2}\geq a_{e-1}\geq ... \geq a_1\geq 0$ (the inequalities being strict if the $a_i$'s are positive), and each $a_i\in O(m)$.
\end{lemma}

\begin{proof}
The $a_i$'s are simply obtained from the $(e-2)$-binomial expansion of $r-m-\binom{m+e-3}{e-1}$ (we can consider them to be all 0's if $r=m-\binom{m+e-3}{e-1}$).
\end{proof}

The following result is due to Stanley, even if its idea was already contained in a paper of Reiten (\cite{Re}).

\begin{lemma}\label{tri}
Given a level algebra with $h$-vector $(1,h_1,...,h_j)$, there exists a Gorenstein algebra (called its {\em trivial extension}) having $h$-vector $H=(1,H_1,...,H_j,H_{j+1})$, where, for each $i=1,2,...j$, we have
$$H_i=h_i+h_{j+1-i}.$$
\end{lemma}

\begin{proof}
See \cite{St}, Example 4.3.
\end{proof}

The following is the main result of this paper.  Notice that once we have fixed the socle degree $e$,  by symmetry it is enough to determine the behavior of the Hilbert function in degrees $i \leq \frac{e}{2}$ as $r \rightarrow \infty$.  Notice also that the following result generalizes Stanley's conjecture when $i = 2$ and $e = 4$, which we proved in \cite{MNZ1}. Also, it greatly generalizes a theorem of Kleinschmidt (\cite{kleinschmidt}, Theorem 1), which supplied a logarithmic estimate for the middle entry, namely:
$$\log f_{e,\left \lfloor \frac{e}{2}\right \rfloor }(r)\sim_r \frac{\left\lfloor \frac{e+1}{2}\right \rfloor }{e-1}\log r.$$
(Recall that two arithmetic functions $f$ and $g$ are asymptotic, i.e., $f(r)\sim_r g(r)$, when $\lim_{r\rightarrow +\infty }\frac{f(r)}{g(r)}=1$. One often simply writes $f(r)\sim_r g(r)$ in place of $f(r)\sim_{r \rightarrow +\infty }g(r)$, since $+\infty $ is the only accumulation point for the natural numbers with respect to the discrete topology they naturally inherit from the reals.) 

A surprising fact is that the asymptotic formula we will  show for
$f_{e,i}(r)$ suddenly increases by a factor of 2 exactly in the
middle (i.e., when $i=\frac{e}{2}$; therefore this pathology occurs
only when the socle degree $e$ is even).

\begin{theorem}  \label{main result} Fix $e$ and $i$.  Then
\[
\displaystyle \lim_{r \rightarrow \infty} \frac{f_{e,i}(r)}{r^{\frac{e-i}{e-1}}} =
\left \{
\begin{array}{ll}
\displaystyle \frac{\left(\left(e-1\right)!\right) ^{\frac{e-i}{e-1}}}{(e-i)!} & \hbox{\hskip 1cm if $i < \frac{e}{2}$} \\ \\
\displaystyle 2 \cdot \frac{\left(\left(e-1\right)!\right) ^{\frac{e/2}{e-1}}}{(e/2)!} & \hbox{\hskip 1cm if $i = \frac{e}{2},$}
\end{array}
\right.
\]
where $f_{e,i}(r)$, as in the above definition, denotes the least possible value that the Hilbert function of a Gorenstein algebra of  codimension $r$ and socle degree $e$ may assume in degree $i$. (Notice that, if $i = \frac{e}{2}$, the  left-hand side of the displayed equation has denominator equal to $r^{\frac{e/2}{e-1}}$.)
\end{theorem}

\begin{proof}
Let $F(r) := f_{e,i}(r)/r^{\frac{e-i}{e-1}}$.  We have to show that the limit exists and is equal to the asserted value.  This was done for $e=4$ and $i=2$ in \cite{MNZ1}, so we will assume that $e \geq 5$.  We will exhibit functions $G$ and $H$ such that, for all $r$, $G(r) \leq F(r) \leq H(r)$ and  both $G$ and $H$ converge to the limit asserted in the theorem.
We begin by producing $G(r)$.

We first assume that $i < \frac{e}{2}$.  Observe that by Theorem \ref{lower bd} (or by Theorem 4 of \cite{MNZ1}) and the fact that $h_1 = r$, we have
\begin{equation} \label{exp of h2}
h_2 \geq \left (r_{(e-1)} \right )^{-1}_{-1} + \left ( r_{(e-1)} \right )^{-(e-2)}_{-(e-3)} \geq \left (r_{(e-1)} \right )^{-1}_{-1} .
\end{equation}
Consider the $(e-2)$-binomial expansion of $h_2$:
\[
(h_2)_{(e-2)} = \binom{\alpha_{e-2}}{e-2} + \binom{\alpha_{e-3}}{e-3} + \dots + \binom{\alpha_1}{1}.
\]
Then again by Theorem \ref{lower bd} we have for $h_3$ that
\[
\begin{array}{rcl}
h_3 & \geq & \displaystyle  \left ( \left(h_2\right )_{(e-2)} \right )^{-1}_{-1} + \left ( \left(h_2\right )_{(e-2)} \right )^{-(e-4)}_{-(e-5)} \\ \\
& \geq & \left ( \left(h_2\right )_{(e-2)} \right )^{-1}_{-1} \\ \\

& \geq & \left ( \left ( \left ( r_{(e-1)} \right )^{-1}_{-1} \right )_{(e-2)} \right )^{-1}_{-1} \hfill \hbox{\hskip 6cm (by (\ref{exp of h2}))} \\ \\
& = & \left ( r_{(e-1)} \right )^{-2}_{-2} \hfill \hbox{\hskip 6cm (by Corollary \ref{minus}).}
\end{array}
\]
Proceeding inductively in the same way, we obtain for $i < \frac{e}{2}$, using Corollary \ref{minus},
\begin{equation} \label{first ineq}
f_{e,i}(r) \geq \left ( r_{(e-1)} \right )^{-(i-1)}_{-(i-1)}.
\end{equation}

Consider the $(e-1)$-binomial expansion of $r$:
\[
r_{(e-1)} = \binom{k}{e-1} + \binom{k_{e-2}}{e-2} + \dots + \binom{k_1}{1}.
\]

Note that $k$ is obtained as a function of $r$.  Thus, invoking (\ref{first ineq}), we obtain
\[
f_{e,i}(r) \geq \binom{k-i+1}{e-i} .
\]

Since $k$ is a function of $r$, and $e$ and $i$ are fixed in advance, $\binom{k-i+1}{e-i}$ is also a function of $r$, which we denote by $G_1(r)$.

Since asymptotically $r \sim_r k^{e-1}/(e-1)!$, we have
\[
k \sim_r r^{\frac{1}{e-1}} \cdot \left(\left(e-1\right)!\right)^{\frac{1}{e-1}},
\]
and so
\[
G_1(r) \sim_r \frac{k^{e-i}}{(e-i)!}  \sim_r \frac{r^{\frac{e-i}{e-1}} \cdot ((e-1)!)^{\frac{e-i}{e-1}}}{(e-i)!}.
\]
Denoting  $G(r) := G_1(r)/r^{\frac{e-i}{e-1}}$,
we see that $G(r) \leq F(r)$ and $G(r)$ has the desired limit when $i < \frac{e}{2}$.

The argument is similar when $i = \frac{e}{2}$, with essentially one difference.  We now have, using Theorem \ref{lower bd}, that
\[
h_{\frac{e}{2}} \geq \left ( \left(h_{\frac{e}{2}-1} \right)_{(e-\frac{e}{2}+1)} \right )^{-1}_{-1} +
\left ( \left(h_{\frac{e}{2}-1}\right )_{(e - \frac{e}{2}+1)} \right )^{-2}_{-1} =
 \left ( \left(h_{\frac{e}{2}-1} \right)_{(\frac{e}{2}+1)} \right )^{-1}_{-1} +
\left ( \left(h_{\frac{e}{2}-1} \right)_{( \frac{e}{2}+1)} \right )^{-2}_{-1} .
\]

Arguing as before, we now obtain
\[
f_{e,\frac{e}{2}} \geq \left ( r_{(e-1)} \right )^{-\frac{e}{2} +1}_{-\frac{e}{2} +1} + \left ( r_{(e-1)} \right )^{-\frac{e}{2} }_{-\frac{e}{2} +1} .
\]

Since asymptotically both terms carry equal weight, we proceed as before with a factor of two, as asserted.

We now want to show the upper bound, by exhibiting a function $H(r)\geq F(r)$ which converges  to the limit of the statement.

Let us write $r$ as in Lemma \ref{r}, and consider the integer $r-m=\binom{m+e-3}{e-1}+\binom{a_{e-2}}{e-2}+\binom{a_{e-3}}{e-3}+...+\binom{a_1}{1}.$ First suppose that $r>m+\binom{m+e-3}{e-1}$, i.e. that $a_{e-2}\geq e-2$.

We construct an $h$-vector $h$ of socle degree $e$ and type $h_{e-1}=r-m$ as follows. For all indices $i$, let $$h_i=\left(\left(h_{e-1}\right)_{(e-1)}\right)^{-(e-1-i)}_{-(e-1-i)}=\binom{m-2+i}{i}+\binom{a_{e-2}-e+i+1}{i-1}+...+\binom{a_{e-i-1}-e+i+1}{0}.$$
In particular, $$h_1=\binom{m-1}{1}+\binom{a_{e-2}-e+2}{0}=(m-1)+1=m.$$

It is easy to see, by the fact that all $a_i$'s are $O(m)$, that $h_{i}\sim_m \frac{m^i}{i!}$.

Furthermore, by Lemma \ref{BG Lemma} and Corollary \ref{minus}, we have that $h_i$ is the minimum possible value of $h$ in degree $i$, given $h_{e-1}$. It is easy to show that this construction guarantees that $h$ be level, since the lex-segment ideal corresponding to $h$ is a level ideal (see, e.g., \cite{BG} or \cite{Za4}).

Hence, by trivial extension (Lemma \ref{tri}), we can construct a Gorenstein $h$-vector $(1,H_1,...,H_e)$ of socle degree $e$, where $H_i=h_i+h_{e-i}$.

In particular, $H_1=h_1+h_{e-1}=m+(r-m)=r$. Also, for all indices $i\leq \frac{e}{2}$, we have
$$H_i\sim_m \frac{m^i}{i!}+\frac{m^{e-i}}{(e-i)!},$$ which is asymptotic to $\frac{m^{e-i}}{(e-i)!}$ if $i< \frac{e}{2}$, and to $2\frac{m^{e/2}}{(e/2)!}$ if $i=\frac{e}{2}$.

Since $m$ is a function of $r$, $H_i$ is also a function of $r$. Also, notice that, asymptotically, $r\sim_r \frac{m^{e-1}}{(e-1)!}$, whence $m\sim_r ((e-1)!)^{\frac{1}{e-1}}r^{\frac{1}{e-1}}$.

Thus, since by definition, $f_{e,i}(r)\leq H_i(r)$, we have
$$\frac{f_{e,i}(r)}{r^{\frac{e-i}{e-1}}}\leq \frac{H_i(r)}{r^{\frac{e-i}{e-1}}},$$
and it is easy to check that the right hand side converges to the desired value for all $i\leq \frac{e}{2}$.

It remains to prove the upper bound when $r$ is of the  form
$r=m+\binom{m+e-3}{e-1}$.

We proceed exactly as before, by starting with a  level $h$-vector
of type $h_{e-1}=r-m=\binom{m+e-3}{e-1}$, and obtaining, by trivial
extensions, a Gorenstein $h$-vector $(1,H_1,...,H_e)$, where $H_i=
\binom{m+i-2}{i}+\binom{m+e-i-2}{e-i}$ if $0 \leq i \leq e-1$. The
only difference is that now
$H_1=\binom{m-1}{1}+\binom{m+e-3}{e-1}=(m-1)+(r-m)=r-1$.

But it is easy to show that  if
$(1,H_1,H_2,...,H_{e-1},1)$ is a Gorenstein $h$-vector, then also $(1,H_1+1,H_2+1,...,H_{e-1}+1,1)$ is always a Gorenstein $h$-vector (for instance using Macaulay's inverse systems; see, e.g., the proof of Proposition 8 in \cite{MNZ1}).

Hence, we have constructed a Gorenstein $h$-vector  of codimension
$r$ also when $r=m+\binom{m+e-3}{e-1}$, and, employing the same
argument as above, we obtain that asymptotically its entries again
satisfy the estimate of the statement, since adding 1 clearly does
not change their asymptotic value.

The proof of the theorem is complete.
\end{proof}

We illustrate the quality of our bounds by an example in which we focus on degrees two and three.

\begin{example} \label{relate results}
  Consider the degrees 2 and 3 entries of a Gorenstein $h$-vector  \linebreak $(1, h_1,
h_2,h_3,\ldots,h_e)$, where $r = h_1 = \binom{m+e-3}{e-1} + m$ for
some integer $m$ satisfying $1 \leq m  \leq e-2$.     Assume that $e \geq 6$.  Note that the construction given in the last part of the proof of Theorem \ref{main result} gives a Gorenstein algebra with $h$-vector
{\small 
\begin{equation} \label{hvtr}
\left ( 1 , \binom{m+e-3}{e-1} +m, \binom{m+e-4}{e-2} + \binom{m}{2} + 1, \binom{m+e-5}{e-3} + \binom{m+1}{3} + 1, \dots \right ) 
\end{equation}
}
One quickly checks that 
\[
h_1 = r = \binom{m+e-3}{e-1} + m = \binom{m+e-3}{e-1} + \binom{e-2}{e-2} + \dots + \binom{e-m-1}{e-m-1} 
\]
so applying Theorem \ref{lower bd} to get a bound for $h_2$, we obtain
\[
\begin{array}{rcl}
\displaystyle h_2 &  \displaystyle \geq & \displaystyle  \left [ \binom{m+e-4}{e-2} + \binom{e-3}{e-3} + \dots + \binom{e-m-2}{e-m-2} \right ] + \left [ \binom{m+e-3-e+2}{e-1-e+3} + 0 \right ] \\
\displaystyle & = &\displaystyle \binom{m+e-4}{e-2} + m + \binom{m-1}{2}
\end{array}
\] 
Since $\binom{m-1}{2} + m = \binom{m}{2} + 1$, we see that for this class of examples the bound for $h_2$ given in Theorem \ref{lower bd} is sharp!

Similarly, let us consider the bound that we obtain for $h_3$. We have already computed in (\ref{hvtr}) the value of $h_3$ obtained in the construction of Theorem \ref{main result}.  To apply Theorem \ref{lower bd}, we need to write the $(e-2)$-binomial expansion of $h_2$.  To that end, suppose that $a \geq 1$ and $k \leq e-2$ are integers satisfying
\[
(e-2) + (e-3) + \dots + (e-k) + a = \binom{m}{2} + 1.
\]
Notice that, since $m\leq e-2$ and $e\geq 6$, such integers $a$ and $k$ always exist.  Hence
\[
h_2 \geq \binom{m+e-4}{e-2} + \binom{e-2}{e-3} + \binom{e-3}{e-4} + \dots + \binom{e-k}{e-k-1} + a
\]
where here we are thinking of $a$ as a sum of binomial coefficients of the form $\binom{e-j}{e-j}$.  Then Theorem \ref{lower bd} gives
\[
\begin{array}{rcl}
h_3 & \geq &  \left [ \binom{m+e-5}{e-3} + \binom{e-3}{e-4} + \dots + \binom{e-k-1}{e-k-2} + a \right ] + \left [ \binom{m+e-4 - (e-4)}{(e-2)-(e-5)} + (k-1) \right ] \\
& = & \left [ \binom{m+e-5}{e-3} + (e-3) + \dots + (e-k-1) + a \right ] + \left [ \binom{m}{3} + (k-1) \right ] \\
& = & \binom{m+e-5}{e-3} + \binom{m}{2} + 1 + \binom{m}{3} \\
& = & \binom{m+e-5}{e-3} + \binom{m+1}{3} + 1.
\end{array}
\]
Hence the bound of Theorem \ref{lower bd} is attained.
Choosing $e = 8$ and $m = 5$ we obtain the example given in the introduction.
\end{example}
{\ }\\
{\bf Acknowledgements.} We thank the anonymous referee for very helpful comments.


\begin{thebibliography}{ll}

\bibitem{bass} H.\ Bass, {\em On the ubiquity of Gorenstein rings},  Math.\ Z.\  {\bf 82}  (1963), 8--28.

\bibitem{BI} D. Bernstein and A. Iarrobino: {\em A non-unimodal graded Gorenstein Artin algebra in codimension five}, Comm.\ in Algebra {\bf 20} (1992), No. 8, 2323-2336.

\bibitem{BG} A.M.\ Bigatti and A.V.\ Geramita: {\em Level Algebras, Lex Segments and Minimal Hilbert Functions}, Comm.\ in Algebra {\bf 31} (2003), 1427-1451.

\bibitem{Bo} M.\ Boij: {\em Graded Gorenstein Artin algebras whose Hilbert functions have a large number of valleys}, Comm.\ in Algebra {\bf 23} (1995), No. 1, 97-103.

\bibitem{BL} M.\ Boij and D.\ Laksov: {\em Nonunimodality of graded Gorenstein Artin algebras}, Proc.\ Amer.\ Math.\ Soc.\ {\bf 120} (1994), 1083-1092.

\bibitem{BH} W.\ Bruns and J.\ Herzog: {\em Cohen-Macaulay rings}, Cambridge studies in advanced mathematics, No.\ 39, Revised edition (1998), Cambridge, U.K..

\bibitem{EI} J.\ Emsalem and A.\ Iarrobino: {\em Some zero-dimensional generic singularities; finite algebras having small tangent space}, Compositio Math.\ {\bf 36} (1978), 145-188.

\bibitem{FL} R.\ Fr\"oberg and D.\ Laksov: {\em Compressed Algebras}, Conference on Complete Intersections in Acireale, Lecture Notes in Mathematics, No. 1092 (1984), 121-151, Springer-Verlag.

\bibitem{Gr} M.\ Green: {\it Restrictions of linear series to hyperplanes, and some results of Macaulay and Gotzmann}, Algebraic curves and projective geometry (1988), 76-86, Trento; Lecture Notes in Math. {\bf 1389} (1989), Springer, Berlin.

\bibitem{Ia1} A. Iarrobino: {\em Compressed Algebras: Artin algebras having given socle degrees and maximal length}, Trans. Amer. Math. Soc. {\bf 285} (1984), 337-378.

\bibitem{IS} A.\ Iarrobino and H.\ Srinivasan: {\em Some Gorenstein Artin algebras of embedding dimension four, I: components of $PGOR(H)$ for $H=(1,4,7,...,1)$}, J.\ of Pure and Applied Algebra {\bf 201} (2005), 62-96.

\bibitem{kleinschmidt}  P.\ Kleinschmidt: {\em \"Uber Hilbert-Funktionen graduierter Gorenstein-Algebren}, Arch.\ Math. {\bf 43} (1984), 501-506.

\bibitem{Ma} F.H.S. Macaulay: {\em The Algebraic Theory of Modular Systems}, Cambridge Univ. Press, Cambridge, U.K. (1916).

\bibitem{MNZ1}
J.\ Migliore, U.\ Nagel and F.\ Zanello:   {\em On the degree two entry of a Gorenstein $h$-vector and a conjecture of
Stanley}, Proc.\ Amer.\ Math.\ Soc. {\bf 136} (2008), No. 8, 2755-2762.

\bibitem{MNZ2}
J.\ Migliore, U.\ Nagel and F.\ Zanello:  {\em A characterization of
Gorenstein Hilbert functions in codimension four with small initial
degree}, Math.\ Res.\ Lett. {\bf 15} (2008), No. 2, 331-349.

\bibitem{Re} I.\ Reiten: {\em The converse to a  theorem of Sharp on Gorenstein modules}, Proc.\ Amer.\ Math.\ Soc.\ {\bf 32} (1972), 417-420.

\bibitem{St} R.\ Stanley: {\em Hilbert functions of graded algebras}, Adv.\ Math.\ {\bf 28} (1978), 57-83.

\bibitem{St2} R.\ Stanley: {\em Combinatorics and Commutative Algebra}, Second Ed., Progress in Mathematics {\bf 41} (1996), Birkh\"auser, Boston, MA.

\bibitem{Za1} F.\ Zanello: {\em Extending the idea of compressed algebra to arbitrary socle-vectors}, J.\ of Algebra {\bf 270} (2003), No. 1, 181-198.

\bibitem{Za2} F.\ Zanello: {\em Extending the idea of compressed algebra to arbitrary socle-vectors, II: cases of non-existence}, J.\ of Algebra {\bf 275} (2004), No. 2, 730-748.

\bibitem{Za3} F. Zanello: {\em Stanley's theorem on codimension 3 Gorenstein $h$-vectors}, Proc.\ Amer.\ Math.\ Soc.\ {\bf 134} (2006), No. 1, 5-8.

\bibitem{Za4} F. Zanello: {\em When is there a unique socle-vector associated to a given $h$-vector?}, Comm. in Algebra {\bf 34} (2006), No. 5, 1847-1860.

\end{thebibliography}
\end{document}